# SHAPE RECONSTRUCTION IN SCATTERING MEDIA WITH VOIDS USING A TRANSPORT MODEL AND LEVEL SETS

OLIVER DORN

ABSTRACT. A two-step shape reconstruction method for diffuse optical tomography (DOT) is presented which uses adjoint fields and level sets. The propagation of near-infrared photons in tissue is modeled by the time-dependent linear transport equation, of which the absorption parameter has to be reconstructed from boundary measurements. In the shape reconstruction approach, it is assumed that the inhomogeneous background absorption parameter and the values inside the obstacles (which typically have a high contrast to the background) are known, but that the number, sizes, shapes, and locations of these obstacles have to be reconstructed from the data. An additional difficulty arises due to the presence of so-called clear regions in the medium. The first step of the reconstruction scheme is a transport-backtransport (TBT) method which provides us with a low-contrast approximation to the sought objects. The second step uses this result as an initial guess for solving the shape reconstruction problem. A key point in this second step is the fusion of the 'level set technique' for representing the shapes of the reconstructed obstacles, and an 'adjoint-field technique' for solving the nonlinear inverse problem. Numerical experiments are presented which show that this novel method is able to recover one or more objects very fast and with good accuracy.

1. INTRODUCTION

1.1. **A novel two-step shape reconstruction scheme.** In a recent paper [12, 13] a two-step shape reconstruction method was introduced for the situation of 2D cross-borehole electromagnetic (EM) imaging. In that paper, electromagnetic cross-borehole tomography was modeled as an inverse problem for the 2D Helmholtz equation. By reformulating the inverse problem as a shape reconstruction problem, and applying a novel level set based shape reconstruction scheme, it was possible to recover nontrivial shapes from noisy, limited-view, multifrequency data.

The method presented in [12, 13] is quite general. In the present paper, we demonstrate how this novel two-step shape reconstruction scheme can be applied to a completely different imaging situation. We are now interested in the imaging of absorption characteristics in human tissue from near-infrared light data. This imaging technique has become known as Diffuse Optical Tomography (DOT), see for example [1]. In contrast to the majority of approaches in DOT which use a diffusion approximation to the linear transport equation for the inversion task, we consider here the (time-dependent) linear transport equation itself, which is generally accepted as being the correct model for the photon propagation in tissue.

This article was written during the stay of the author at the Mathematical Sciences Research Institute (MSRI), Berkeley, from October to December 2001. The research presented here is supported in part by the NSERC CRD Grant 80357, and research during the stay at MSRI is supported by NSF grant DMS-9810361.

The inversion method is a two-step shape reconstruction scheme which tries to recover the unknown shapes of abnormalities in the absorption coefficient in the tissue. The first step consists of an inverse scattering algorithm for DOT which is called a 'Transport-Backtransport' (TBT) method and which was originally introduced in [9, 10]. In the situation of DOT, this TBT algorithm yields already during the early iterations a relatively good approximation to the locations of the unknown obstacles, but requires many more iterations in order to get better estimates for the correct extension and shapes of these objects. Therefore, the result of these early sweeps is only used as an initial guess for initiating the second step of the two-step scheme, which is a level set based shape reconstruction method. This second step calculates successive corrections for the shapes of the obstacles which aim to reduce the misfit in the data. Here, the level set representation of the shapes is chosen because it enables us to easily describe and keep track of complicated geometries which might arise during the evolution process.

Both, the initialization step as well as the level set based shape evolution step employ an 'adjoint-field technique' for the inversion which has the very useful property that the inverse problem can be solved approximately by making two uses of the same forward modeling code. Using a somewhat oversimplified description of our technique, the updates to the level set function are obtained by first making one pass through the code using the parameter distribution (which is here the absorption parameter) corresponding to the latest best guess of the level set function, and then another pass with the adjoint operator applied to the differences in computed and measured data. Then the results of these two calculations are combined to determine updates to the level set function. The resulting procedure is iterative, and can be applied successively to parts of the data, e.g., data associated with one source location can be used to update the model before other source locations are considered. Similar techniques have been successfully applied recently as part of iterative nonlinear inversion schemes in [6, 10, 11, 32, 44]. For alternative inversion schemes in DOT see for example [1, 17, 18, 20, 25, 27, 30, 45] and the references therein.

1.2. **The level set method.** The *level set method* was originally developed by Osher and Sethian for describing the motion of curves and surfaces [36, 41]. Since then, it has found applications in a variety of quite different situations. Examples are image enhancement, computer vision, interface problems, crystal growth, or etching and deposition in the microchip fabrication. For an overview see [42].

The idea of using a level set representation as part of a solution scheme for inverse problems involving obstacles was first suggested by Santosa in [40]. In that paper, two linear inverse problems, a deconvolution problem and the problem of reconstructing a diffraction screen, are solved by employing an optimization approach as well as a time evolution approach using the level set technique for describing the shapes. Santosa also outlines in that paper how the two presented methods can be generalized to nonlinear shape recovery problems. More recently, a related method was applied to a nonlinear shape recovery problem by Litman *et al.* in [28]. In that work, an inverse transmission problem in free space is solved by a controled evolution of a level set function. This evolution is governed by a Hamilton-Jacobi type equation, whose velocity function has to be determined properly in order to minimize a given cost functional.

In Dorn *et al.* [12, 13] a two-step shape reconstruction algorithm was presented for the situation of electromagnetic cross-borehole tomography in 2D. In contrary to the work mentioned above, this approach is not based on a Hamilton-Jacobi type equation, but uses an optimization approach instead. Moreover, a very specific inversion routine, an adjoint-field technique, is employed which calculates an update for the most recent level set function directly by just solving one forward and one adjoint problem on the computer. The present paper is based on this approach, and a detailed discussion will be given in the following sections.

As further examples for related and very interesting recent approaches we also want to mention here the work described in Burger [2], Elangovan et al [16], Hintermüller et al [21], Ito et al [22], Jackowska-Strumillo et al [23], Osher et al [37], Ramananjaona et al. [38], Sethian et al [43], and Zhao et al [46].

The paper is organized as follows. In section 2 we will present the basic equations of DOT in a form convenient for development of the shape reconstruction technique. In section 3 we formulate the shape reconstruction problem and introduce the level set formulation of this problem. In section 4, we derive the basic shape reconstruction algorithm using level sets and adjoint fields. Section 5 describes shortly the procedure employed here for finding a good starting guess for our shape recontruction algorithm. In Section 6 numerical experiments are presented which demonstrate the performance of the algorithm in different situations. The final section summarizes the results of this paper and indicates some directions for future research.

## 2. The physical experiment in DOT

### 2.1. The linear transport equation.
We model the propagation of photons in the medium by the time-dependent linear transport equation

$$\frac{\partial u}{\partial t} + \theta \cdot \nabla u(x,\theta,t) + (a(x)+b(x))u(x,\theta,t) - b(x)\int_{S^{n-1}} \eta(\theta \cdot \theta')u(x,\theta',t)d\theta'$$

$$= q(x,\theta,t) \quad \text{in} \quad \Omega \times S^{n-1} \times [0,T] \tag{1}$$

with initial condition

$$u(x,\theta,0) = 0 \quad \text{in} \quad \Omega \times S^{n-1} \tag{2}$$

and boundary condition

$$u(x,\theta,t) = 0 \quad \text{on} \quad \Gamma_-. \tag{3}$$

Here,

$$\Gamma_\pm := \left\{(x,\theta,t) \in \partial\Omega \times S^{n-1} \times [0,T], \quad \pm\nu(x) \cdot \theta > 0\right\}.$$

$\Omega$ is a convex, compact domain in $\mathbb{R}^n$, $n=2,3$, with smooth boundary $\partial\Omega$. In our numerical experiments, we will only consider the case $n=2$, but the algorithm extends in a straightforward way to $n=3$. $\nu(x)$ denotes the outward unit normal to $\partial\Omega$ at the point $x \in \partial\Omega$, and $u(x,\theta,t)$ describes the density of particles (photons) which travel in $\Omega$ at time $t$ through the point $x$ in the direction $\theta$. The velocity $c$ of the particles is assumed to be normalized to $c = 1\text{cms}^{-1}$ and has been neglected in the formulation of (1).

$a(x)$ is the absorption cross-section (in short 'absorption'), $b(x)$ is the scattering cross-section, and $\mu(x) := a(x)+b(x)$ is the total cross-section or attenuation. These parameters are assumed to be real, strictly positive functions of the position $x$. The quantity $\mu^{-1}$ is the mean free path of the photons. Typical values in DOT are $a \approx 0.1 - 1.0 \text{cm}^{-1}$, $b \approx 100.0 - 200.0 \text{cm}^{-1}$, $\mu^{-1} \approx 0.005 - 0.01 \text{cm}$, see for example [24]. The scattering function $\eta(\theta \cdot \theta')$ describes the probability for a particle entering a scattering process with the direction of propagation $\theta'$ to leave this process with the direction $\theta$. It is normalized to

$$\int_{S^{n-1}} \eta(\theta \cdot \theta') d\theta = 1, \tag{4}$$

which expresses particle conservation in pure scattering events. The dot-product in the argument indicates that $\eta$ depends only on the cosine of the scattering angle $\cos \vartheta = \theta \cdot \theta'$ and, in particular, is independent of the location of the scattering event $x$. These latter assumptions simplify the following calculations, but the reconstruction method can be extended to more general scattering functions as well. In our numerical experiments, we will use the following Henyey-Greenstein scattering function [24]

$$\eta(\theta \cdot \theta') = \frac{1 - g^2}{2(1 + g^2 - 2g \cos \vartheta)^{3/2}}, \tag{5}$$

with $-1 < g < 1$. The parameter $g$ in (5) is the mean cosine of the scattering function. Values of $g$ close to one indicate that the scattering is primarily forward directed, whereas values close to zero indicate that scattering is almost isotropic. In our numerical experiments, we will choose $g$ to be 0.9, which is a typical value for DOT.

The initial condition (2) indicates that there are no photons moving inside of $\Omega$ at the starting time of our experiment. The boundary condition (3) indicates that during the experiment no photons enter the domain $\Omega$ from the outside. All photons inside of $\Omega$ originate from the source $q$ which, however, can be situated at the boundary $\partial \Omega$. We mention that an alternative way of describing laser sources at the boundary is to use an inhomogeneous boundary condition instead of (3) and a zero source $q = 0$. Both choices are equivalent, and the derivation of the TBT-method as well as the shape reconstruction method are almost identical for these two possibilities.

We consider the problem (1)-(3) for $p$ different sources $q_j$, $j = 1, \ldots, p$, positioned at $\partial \Omega$. Typical sources in applications are delta-like pulses transmitted at time $t = 0$ at the position $s_j \in \partial \Omega$ into the direction $\theta_j$, which can be described by the distributional expressions

$$\tilde{q}_j(x, \theta, t) = \delta_x(s_j) \delta_t(0) \delta_\theta(\theta_j), \qquad j = 1, \ldots, p, \tag{6}$$

with $\nu(s_j) \cdot \theta_j < 0$. We assume that a given source $q_j$ gives rise to the physical fields $\tilde{u}_j(x, \theta, t)$ which are solutions of (1)-(3). Our measurements consist of the outgoing flux across the boundary $\partial \Omega$ which has with (3) the form

$$\tilde{G}_j(x, t) = \int_{S^{n-1}} \nu(x) \cdot \theta \tilde{u}_j(x, \theta, t) d\theta \quad \text{on} \quad \partial \Omega \times [0, T], \tag{7}$$

for $j = 1, \ldots, p$. With these definitions, we are now able to formulate the inverse problem which we will consider in the following.

**Definition 2.1.** Inverse problem of DOT: *Assume that we measure for $p$ different sources $q_j$, $j = 1, \ldots, p$, the data (7). Provided with this information and knowing the scattering functions $\eta$ and $b$, how can we reconstruct the coefficient $a(x)$ of (1) inside of $\Omega$?*

In general, also the scattering parameter $b(x)$ is part of the inversion problem. However, for simplicity, we will restrict ourselves in this paper to the determination of only one parameter, namely $a(x)$, and will assume that the scattering parameter $b(x)$ is given. Extensions of this method will be addressed in our future research.

2.2. **Clear regions.** An additional difficulty arises in our application of DOT due to the presence of so-called 'clear regions' in the imaging domain. These are regions where the scattering coefficient $b(x)$ is very small or zero, such that the dominant behavior of the traveling photons inside these regions is transport rather than diffusion. Most of the imaging methods in DOT are purely based on the diffusion approximation to the linear transport equation, such that they cannot be used in applications where clear regions are present. Clear regions play an important role for example when imaging the head of neonates or in functional brain imaging. See [1, 26, 35, 39].

Although we assume in this paper that the locations of the clear regions are known a priori (which in realistic applications is not always the case such that the shapes of the clear regions become part of the imaging task themselves), the presence of these regions complicates the numerical inversion process.

One problem caused by the clear layers is their effect on the sensitivity structure of the data with respect to the unknown parameters. The clear regions (and, in particular, the clear layers which are surrounding the region of interest, see figure 1) act like 'wave guides' for the photons, since the propagation is much faster inside these layers than it is in the diffusive regions. Therefore, the photons which have travelled through these clear layers dominate in particular the early parts of the time-dependent measurements. Unfortunately, these photons do not carry much information related to the internal structure of the medium. In absence of clear layers, the early photons are usually assumed to yield the reconstructions with the best resolution. If clear layers are present, however, these early photons are now mixed with those photons that travelled through the clear layers, such that it is difficult to extract the necessary information about the interior domain from the early data. Therefore, we will in our numerical reconstructions only use photons which arrived at later times at the receiver locations. Our sensitivity analysis indicates that these later photons are less affected by the clear layers. We will show that the information of the later arriving photons is sufficient for achieving good reconstructions of the absorption distribution inside the regions of interest, which might however have a slightly reduced resolution. See also the more detailed discussion of sensitivity functions in section 4.7.

## 3. The shape reconstruction problem

In this section we formulate the shape reconstruction problem which we want to solve, and cast it in a form which makes use of the level set representation of the domains.

3.1. **The shape reconstruction problem in DOT.** To start with we introduce some terminology which we will use throughout the paper.

**Definition 3.1.** *Let us assume that we are given a constant $\hat{a} > 0$, and a bounded function $a_b : \Omega \to \mathbf{R}$. We call a pair $(D, a)$, which consists of a compact domain*

$D \subset\subset \Omega$ and a bounded function $a : \Omega \to \mathbb{R}$, admissible *if we have*

$$a|_D = \hat{a}, \quad a|_{\Omega \setminus D} = a_{\mathrm{b}}|_{\Omega \setminus D}. \tag{8}$$

In other words, a pair $(D, a)$ is *admissible* if $a$ is equal to a preassigned constant value $\hat{a}$ inside of $D$, and equal to the preassigned *background absorption* $a_{\mathrm{b}}$ outside of $D$. The domain $D$ is called the *scattering domain*.

**Remark 3.1.** *For an admissible pair $(D, a)$, and for given $\hat{a}$, $a_{\mathrm{b}}$, the absorption $a$ is uniquely determined by $D$.*

With this definition, we can now formulate the *shape reconstruction problem*.

*Shape reconstruction problem.* Let us assume that we are given a constant $\hat{a} > 0$, a bounded function $a_{\mathrm{b}} : \Omega \to \mathbb{R}$, and some data $\tilde{G}_j$ as in (7). Find a domain $\tilde{D}$ such that the admissible pair $(\tilde{D}, \tilde{a})$ reproduces the data, i.e. (7) holds with $\tilde{u}_j$ being a solution of (1)-(3) with $a = \tilde{a}$ for $j = 1, \ldots, p$.

Using the same notation and assumptions as in definition 3.1, we want to formulate another inverse problem which we will call the *inverse scattering problem* and which will play an important part when solving the shape reconstruction problem.

*Inverse Scattering Problem.* Let us assume that we are given a bounded function $a_{\mathrm{b}} : \Omega \to \mathbb{R}$, and some data $\tilde{G}_j$ as in (7). Find a bounded function $\tilde{a}_{\mathrm{s}} : \Omega \to \mathbb{R}$ such that $\tilde{a} = a_{\mathrm{b}} + \tilde{a}_{\mathrm{s}}$ reproduces the data, i.e. (6) holds with $\tilde{u}_j$ given by (1)-(3) and with $a = \tilde{a}$ for $j = 1, \ldots, p$.

The inverse scattering problem gives rise to the decomposition

$$a = a_{\mathrm{b}} + a_{\mathrm{s}} \quad \text{in } \Omega \tag{9}$$

In other words, the absorption distribution $a$ is decomposed into the background distribution $a_{\mathrm{b}}$ and the perturbation $a_{\mathrm{s}}$ which we will refer to in the following as the *obstacle*.

Solving the shape reconstruction problem requires only to find the *shape* of the domain $\tilde{D}$, since the function $\tilde{a}$ is then uniquely determined by (8). Solving the inverse scattering problem, on the other hand, amounts to finding the entire function $\tilde{a}_{\mathrm{s}}$ from the given data, which is much harder to do. However, it will turn out that finding a good *approximate* solution of the inverse scattering problem is much easier and faster to achieve and will provide us with a very good initial guess for starting our shape reconstruction routine.

Definition 3.1 allows us to formulate a first version of the strategy which we want to use for solving the shape reconstruction problem.

*Strategy for solving the shape reconstruction problem.* Construct a series of admissible pairs $(D^{(n)}, a^{(n)})$, $n = 0, 1, 2, \ldots$, such that the misfit between the physically measured data and the calculated data corresponding to $(D^{(n)}, a^{(n)})$ decreases with increasing $n$, and ideally, i.e. in absence of noise, tends to zero in the limit $n \to \infty$. Use the *approximate* solution of the inverse scattering problem to initialize this series by determining a good starting element $(D^{(0)}, a^{(0)})$.

3.2. **The domains** $D^{(n)}$**.** In our numerical examples, each of the domains $D^{(n)}$ which we are looking for can be given as a collection of a finite number $L_n$ of disjoint, compact subdomains $D_l^{(n)}$, $l = 1, \ldots, L_n$, with

$$D^{(n)} = \bigcup_{l=1}^{L_n} D_l^{(n)}, \qquad D_l^{(n)} \cap D_{l'}^{(n)} = \emptyset \quad \text{for } l \neq l'. \tag{10}$$

The shapes of these subdomains $D_l^{(n)}$ can in principle be arbitrary. In particular, they are allowed to be multiply connected, and to enclose some 'cavities' or 'holes' filled with background material, and moreover the number $L_n$ of these subdomains might (and usually does) vary with the iteration number $n$. For the derivation of the inversion method, we assume that the boundaries $\partial D_l^{(n)}$ of these domains are sufficiently smooth.

It is essential for the success and the efficiency of the reconstruction scheme to have a good and flexible way of keeping track of the shape evolution during the reconstruction process. The method we have chosen in our reconstruction algorithm is a level set representation of the shapes as it was suggested by Santosa [40]. This representation has the advantage that the level set functions, which are in principle only used for representing the shapes, can in a natural way be made part of the reconstruction scheme itself. Doing so, it is not necessary anymore to refer to the shapes of the domains until the reconstruction process is completed. The final shape is then recovered from the representing level set function easily. In the following we will discuss in a more formal way how this can be achieved.

3.3. **Level set representation of the domains** $D^{(n)}$**.** Assume that we are given a domain $D \subset\subset \Omega$. The characteristic function $\chi_D : \Omega \to \{0, 1\}$ is defined in the usual way as

$$\chi_D(x) = \begin{cases} 1 & , \quad x \in D \\ 0 & , \quad x \in \Omega \backslash D. \end{cases} \tag{11}$$

**Definition 3.2.** *We call a function $\phi : \Omega \to \mathbb{R}$ a level set representation of $D$ if*

$$\chi_D(x) = \Psi_\phi(x) \qquad \text{on} \quad \Omega \tag{12}$$

*where $\Psi_\phi : \Omega \to \{0, 1\}$ is defined as*

$$\Psi_\phi(x) = \begin{cases} 1 & , \quad \phi(x) \leq 0 \\ 0 & , \quad \phi(x) > 0. \end{cases} \tag{13}$$

For each function $\phi : \Omega \to \mathbb{R}$ there is a domain $D$ associated with $\phi$ by (12),(13) which we call $D[\phi]$. It is clear that different functions $\phi_1, \phi_2$, $\phi_1 \neq \phi_2$, can be associated with the same domain $D[\phi_1] = D[\phi_2]$, but that different domains cannot have the same level set representation. Therefore, we can use the level set representation for unambiguously specifying a domain $D$ by any one of its associated level set functions.

The boundary $\Gamma = \partial D[\phi]$ of a domain $D[\phi]$, represented by the level set function $\phi$, is defined as

$$\begin{aligned} \Gamma = \{x \in \Omega : \quad &\text{for all } \rho > 0 \text{ we can find } x_1, x_2 \in B_\rho(x) \\ &\text{with } \phi(x_1) > 0 \text{ and } \phi(x_2) < 0 \} \end{aligned} \tag{14}$$

**Definition 3.3.** *We call a triple $(D, a, \phi)$, which consists of a domain $D \subset\subset \Omega$ and bounded functions $a, \phi : \Omega \to \mathrm{I\!R}$, admissible if the pair $(D, a)$ is admissible in the sense of definition 3.1, and $\phi$ is a valid level set representation of $D$.*

**Remark 3.2.** *For an admissible triple $(D, a, \phi)$, and for given $\hat{a}$, $a_{\mathrm{b}}$, the pair $(D, a)$ is uniquely determined by $\phi$.*

We use these definitions to reformulate our shape reconstruction problem.

*Level set formulation of the shape reconstruction problem.* Given a constant $\hat{a} > 0$, a background distribution $a_{\mathrm{b}}$, and some data $\tilde{G}$ as in (7). Find a level set function $\tilde{\phi}$ such that the corresponding admissible triple $(\tilde{D}, \tilde{a}, \tilde{\phi})$ reproduces the data, i.e. (7) holds with $\tilde{u}_j$ being a solution of (1)-(3) with $a = \tilde{a}$ for $j = 1, \ldots, p$.

The strategy for solving this shape reconstruction problem has to be reformulated, too. It reads now as follows.

*Strategy for solving the reformulated shape reconstruction problem.* Construct a series of admissible triples $(D^{(n)}, a^{(n)}, \phi^{(n)})$, $n = 0, 1, 2, \ldots$, such that the misfit between the data (7) and the calculated data corresponding to $(D^{(n)}, a^{(n)}, \phi^{(n)})$ decreases with increasing $n$, and ideally, i.e. in absence of noise, tends to zero in the limit $n \to \infty$. For finding this series we only have to keep track of $a^{(n)}$ and $\phi^{(n)}$, but not of $D^{(n)}$. The function $a^{(n)}$ is needed in each step for solving a forward problem (1)-(3), and a corresponding adjoint problem. The knowledge of $\phi^{(n)}$ is used in each step to determine $a^{(n)}$. The final level set function $\phi^{(N)}$, which satisfies some stopping criterion, is used to recover the final shape $D^{(N)}$ via (12).

## 4. Solving the shape reconstruction problem

In this section we derive the basic shape reconstruction method which uses adjoint fields and the level set representation introduced above. The initializing procedure for this reconstruction routine will be discussed in section 5.

### 4.1. Function spaces.
We want to specify now the function spaces which we will be working with. The main objective of this section is to introduce the inner products on these function spaces, which will become important when defining the adjoint linearized operators in the following sections.

The natural spaces in the framework of the linear transport equation are $L^1$-spaces for the fields $u$, and $L^\infty$-spaces for the parameters $a$ and $b$, due to the physical significance of the corresponding norms. This is because the $L^1$-norm of a function $u(x, \theta, t_0)$ describes physically the number of particles which are in the system at a given time $t_0$. However, in our derivation of the shape reconstruction method it will be more convenient to work with $L^2$-spaces instead, since deriving the adjoint operators used in the inversion routine is more straightforward when using these spaces. See for example [5, 7, 9, 10] for more information concerning $L^1$-spaces in linear transport theory. The spaces which we will consider are

$$Y = L^2(\Omega \times S^{n-1} \times [0, T]), \quad U = L^2(\Omega \times S^{n-1} \times [0, T]), \qquad (15)$$

$$Z = L^2(\partial\Omega \times [0, T]), \quad F = L^2(\Omega), \quad \Phi = L^2(\Omega), \qquad (16)$$

equipped with the standard inner products. Here, $Y$ is the space of sources $q$, $U$ is the space of states $u$, $Z$ is the space of measurements or data $G$, $F$ is the space of parameters $a$, and $\Phi$ is the space of level set functions $\phi$. We mention that in fact our level set functions $\phi$ will be assumed to have a certain regularity, since in the derivation of the reconstruction scheme we will make use of first order derivatives of these functions.

4.2. **Operators.** In the following, we will introduce some operators which will enable us to formulate the shape reconstruction problem in a way suitable for deriving the inversion algorithm.

Given a constant $\hat{a}$ and a bounded function $a_{\mathrm{b}} : \Omega \to \mathbb{R}$. Then, with each level set function $\phi \in \Phi$ a uniquely determined 'absorption perturbation' or 'obstacle' $\Lambda(\phi)$ is associated by putting

$$\Lambda(\phi)(x) = \begin{cases} \hat{a} - a_{\mathrm{b}}(x) & , \quad \phi(x) \leq 0 \\ 0 & , \quad \phi(x) > 0. \end{cases} \quad (17)$$

With (13) we can write this also as

$$\Lambda(\phi)(x) = \Psi_\phi(x)(\hat{a} - a_{\mathrm{b}}(x)) \quad , \quad x \in \Omega. \quad (18)$$

Notice that the operator $\Lambda$ is chosen such that the triple $(D, a, \phi)$ with $a = a_{\mathrm{b}} + \Lambda(\phi)$ and domain $D[\phi]$ forms an admissible triple $(D, a, \phi)$ in the sense of definition 3.3. Moreover, for $(D, a, \phi)$ an admissible triple, we see that $\Lambda(\phi)$ is just the absorption perturbation $a_{\mathrm{s}}$

$$\Lambda(\phi)(x) = a_{\mathrm{s}}(x) = \chi_D(x)(\hat{a} - a_{\mathrm{b}}(x)), \quad x \in \Omega. \quad (19)$$

Let us assume now that we are given a background absorption $a_{\mathrm{b}}$ and that we have collected some data $\tilde{G}_j$ which correspond to the 'true' absorption distribution

$$\tilde{a} = a_{\mathrm{b}} + \tilde{a}_{\mathrm{s}}, \quad (20)$$

where $\tilde{a}_{\mathrm{s}}$ is the 'true' obstacle. We consider the following *measurement operators*

$$G_j : F \longrightarrow Z, \quad G_j(a_{\mathrm{s}})(x,t) := \int_{S^{n-1}} \nu(x) \cdot \theta \, u_j(x, \theta, t) \, d\theta, \quad (21)$$

$j = 1, \ldots, p$, where $u_j$ solves (1)-(3) with the source $q_j \in Y$ and the parameters $a = a_{\mathrm{b}} + a_{\mathrm{s}}$. For the 'true' parameters $\tilde{a}_{\mathrm{s}}$ we ask (21) to coincide with the data $\tilde{G}_j$

$$G_j(\tilde{a}_{\mathrm{s}})(x,t) = \tilde{G}_j(x,t) \quad \text{for } j = 1, \ldots, p \quad (22)$$

on $\partial\Omega \times [0, T]$. Assuming $\tilde{G}_j \in Z$ we want to determine $\tilde{a}_{\mathrm{s}}$ such that (22) is valid. For given data $\tilde{G}_j$, $j = 1, \ldots, p$, we define furthermore the *residual operators*

$$R_j : F \longrightarrow Z, \quad R_j(a_{\mathrm{s}})(x,t) = G_j(a_{\mathrm{s}})(x,t) - \tilde{G}_j(x,t). \quad (23)$$

From (22) it follows that for the 'true' obstacle the residuals vanish,

$$R_j(\tilde{a}_{\mathrm{s}}) = 0 \quad \text{for } j = 1, \ldots, p, \quad (24)$$

if the data are noise-free.

The *forward operators* $T_j$ which map a given level set function $\phi \in \Phi$ into the corresponding mismatch in the data are defined by

$$T_j : \Phi \longrightarrow Z_j \quad , \quad T_j(\phi) = R_j(\Lambda(\phi)) \quad (25)$$

for $j = 1, \ldots, p$. The goal is to find a level set function $\tilde{\phi} \in \Phi$ such that

$$T_j(\tilde{\phi}) = 0 \qquad \text{for } j = 1, \ldots, p. \tag{26}$$

We mention that all three operators $\Lambda$, $R_j$ and $T_j$ are nonlinear.

4.3. **Linearized operators.** For the derivation of the shape reconstruction algorithm, we will need expressions for the linearized operators corresponding to the nonlinear operators introduced above, and for their adjoints with respect to the given inner products. In this section, we define the linearized operators, and expressions for their adjoints are derived in section 4.5.

Analogous to Santosa [40] it is shown that, for a homogeneous background $a_{\text{b}}$, the infinitesimal response $\delta a(x)$ in the obstacle $a_{\text{s}}(x)$ to an infinitesimal change $\delta\phi(x)$ of the level set function $\phi(x)$ has the form

$$\delta a(x) = - [\hat{a} - a_{\text{b}}] \frac{\delta\phi(x)}{|\nabla\phi(x)|}\bigg|_{x \in \partial D[\phi]}. \tag{27}$$

The function $\delta a$ in (27) can be interpreted as a 'surface measure' on the boundary $\Gamma = \partial D[\phi]$. Similar to (27), we define the linearized operator $\tilde{\Lambda}'[\phi]$ by

$$\left(\tilde{\Lambda}'[\phi]\delta\phi\right)(x) = - [\hat{a} - a_{\text{b}}(x)] \frac{\delta\phi(x)}{|\nabla\phi(x)|} \hat{\delta}_\Gamma(x) \tag{28}$$

where $\hat{\delta}_\Gamma(x)$ denotes the Dirac delta distribution concentrated on $\Gamma = \partial\Omega[\phi]$. In this interpretation, (27) describes an infinitesimal 'surface load' of absorption on $\Gamma$ which has to be recovered from the mismatch in the data.

However, the expression on the right hand side of (28) is not an element of $F$ which causes problems when we want to calculate the inner products defined in section 4.1. Therefore, we will replace the operator (28) by an approximation which maps from $\Phi$ into $F$ and which will be more convenient for the derivation of the reconstruction method.

For a given level set function $\phi \in \Phi$, let $\Gamma = \partial D[\phi]$ and $B_\rho(\Gamma) = \cup_{y \in \Gamma} B_\rho(y)$ a small finite width neighborhood of $\Gamma$ with some given constant $0 < \rho \ll 1$. The (approximated) linearized operator $\Lambda'[\phi]$ is defined as

$$\Lambda'[\phi] : \Phi \longrightarrow F, \quad (\Lambda'[\phi]\delta\phi)(x) = - [\hat{a} - a_{\text{b}}(x)] \frac{\delta\phi(x)}{|\nabla\phi(x)|} C_\rho(\Gamma) \chi_{B_\rho(\Gamma)}(x) \tag{29}$$

with $\chi_{B_\rho(\Gamma)}(x)$ given by (11). Here, $C_\rho(\Gamma) = L(\Gamma)/\text{Vol}(B_\rho(\Gamma))$ where $L(\Gamma) = \int_\Omega \hat{\delta}_\Gamma(x)dx$ is the length of the boundary $\Gamma$, and $\text{Vol}(B_\rho(\Gamma)) = \int_\Omega \chi_{B_\rho(\Gamma)}(x)dx$ is the volume of $B_\rho(\Gamma)$. For a very small $\rho$ we will get a very large weight $C_\rho(\Gamma)$, whereas for increasing $\rho$ this weight $C_\rho(\Gamma)$ decreases accordingly. The operator defined in (29) maps now from $\Phi$ into $F$ such that we can make use of the inner products defined on these spaces.

We already mentioned that the term $|\nabla\phi(x)|$ in (27), (28), as well as the derivation of these expressions, implies some regularity constraint on $\phi$. For example, $\phi \in C^1$ would be possible. Another possibility would be to use a 'signed distance function' as a standard representation of the boundary [42]. We do not want to specify the regularity of $\phi$ at this point, but assume instead that it is 'sufficiently smooth' for our purposes.

The linearized residual operator $R'_j[a_{\rm s}]$ is defined by

$$R'_j[a_{\rm s}] \,:\, F \longrightarrow Z, \qquad R'_j[a_{\rm s}]\,\delta a \;=\; \int_{S^{n-1}} \nu \cdot \theta\, v_j(x,\theta,t)\, d\theta \tag{30}$$

where $v_j$ solves the problem

$$\frac{\partial v_j}{\partial t} + \theta \cdot \nabla v_j(x,\theta,t) + (a(x)+b(x))\,v_j(x,\theta,t) - b(x) \int_{S^{n-1}} \eta(\theta \cdot \theta')v_j(x,\theta',t)d\theta'$$

$$= -\delta a(x)\, u_j(x,\theta,t) \quad \text{in} \quad \Omega \times S^{n-1} \times [0,T] \tag{31}$$

with initial condition

$$v_j(x,\theta,0) = 0 \quad \text{in} \quad \Omega \times S^{n-1} \tag{32}$$

and boundary condition

$$v_j(x,\theta,t) = 0 \quad \text{on} \quad \Gamma_-. \tag{33}$$

Here, $u_j$ is the solution of (1)–(3) with source $q_j$ and with $a = a_{\rm b} + a_{\rm s}$. Note that (31) has the same structure as (1). The 'primary sources' $q$ of (1) are now replaced by 'secondary sources' which are caused by the parameter perturbation $\delta a$. The representation (30) can be derived by perturbing

$$a_{\rm s} \to a_{\rm s} + \delta a \quad,\quad u_j \to u_j + v_j, \tag{34}$$

plugging this into (1) and neglecting terms which are of higher than linear order in the perturbations $\delta a$, $v_j$. A rigorous derivation of (30) in the framework of $L^1$-spaces can be found in [7].

As our third linearized operator, we introduce the linearized forward operator $T'_j[\phi]$ by putting

$$T'_j[\phi] : \Phi \longrightarrow Z_j \quad,\quad T'_j[\phi]\delta\phi \;=\; R'_j[\Lambda(\phi)]\,\Lambda'[\phi]\delta\phi. \tag{35}$$

All three operators $\Lambda'[\phi]$, $R'_j[a_{\rm s}]$, and $T'_j[\phi]$ are linear.

4.4. **A nonlinear Kaczmarz-type approach.** The algorithm works in a 'single-step fashion' as follows. Instead of using the data (22) for all sources simultaneously, we only use the data for one source at a time while updating the linearized residual operator after each determination of the corresponding incremental correction $\delta\phi$.

To be more specific, let us assume that we are given a level set function $\phi^{(n)}(x)$ and an obstacle $a_{\rm s}^{(n)}(x)$ such that $(D^{(n)}, a_{\rm b} + a_{\rm s}^{(n)}, \phi^{(n)})$ forms an admissible triple in the sense of definition 3.3. Using a data set $\tilde{G}_j$ corresponding to the fixed source position $q_j$, we want to find an update $\delta\phi^{(n)}$ to $\phi^{(n)}$ such that for the admissible triple

$$\left(D^{(n+1)}, a_{\rm b} + a_{\rm s}^{(n+1)}, \phi^{(n+1)}\right) := \tag{36}$$

$$\left(D[\phi^{(n)} + \delta\phi^{(n)}], a_{\rm b} + \Lambda(\phi^{(n)} + \delta\phi^{(n)}), \phi^{(n)} + \delta\phi^{(n)}\right)$$

the residuals in the data corresponding to this source vanish

$$T_j(\phi^{(n+1)}) \;=\; T_j(\phi^{(n)} + \delta\phi^{(n)}) \;=\; 0. \tag{37}$$

Applying a Newton-type approach, we get from (37) a correction $\delta\phi^{(n)}$ for $\phi^{(n)}$ by solving

$$T'_j[\phi^{(n)}]\delta\phi^{(n)} \;=\; -T_j(\phi^{(n)}) \;=\; -\left(M_j u_j - \tilde{G}_j\right) \tag{38}$$

where $u_j$ satisfies (1)–(3) with

$$a(x) = a_{\rm b}(x) + \Lambda(\phi^{(n)})(x) = \begin{cases} \hat{a} & , \quad x \in D[\phi^{(n)}] \\ a_{\rm b}(x) & , \quad x \in \Omega \setminus D[\phi^{(n)}]. \end{cases} \tag{39}$$

Since we have only few data given for one source and one frequency, equation (38) usually will have many solutions (in the absence of noise), such that we have to pick one according to some criterion. We choose to take that solution which minimizes the $L^2$-norm of $\delta\phi^{(n)}$

$$\text{Min } \|\delta\phi^{(n)}\|_2 \quad \text{subject to} \quad T'_j(\phi^{(n)})\delta\phi^{(n)} = -\left(M_j u_j - \tilde{G}_j\right). \tag{40}$$

This solution can be formulated explicitly. It is

$$\delta\phi^{(n)}_{\rm MN} = -T'_j[\phi^{(n)}]^* \left(T'_j[\phi^{(n)}]T'_j[\phi^{(n)}]^*\right)^{-1} \left(M_j u_j - \tilde{G}_j\right), \tag{41}$$

where $T'_j[\phi^{(n)}]^*$ denotes the adjoint operator to $T'_j[\phi^{(n)}]$.

After correcting $\phi$ by $\phi \to \phi + \delta\phi_j$, where $\delta\phi_j$ is given by (41), we use the updated residual equation (38) to compute the next correction $\delta\phi_{j'}$. Doing this for one equation after the other, until each of the sources $q_j$ has been considered exactly once, will yield one complete sweep of the algorithm. This procedure is similar to the Kaczmarz method for solving linear systems, or the algebraic reconstruction technique (ART) in x-ray tomography [31] and the simultaneous iterative reconstruction technique (SIRT) as presented in [8]. For related approaches in a variety of imaging problems we refer to [6, 10, 11, 32, 33, 34, 44].

4.5. **The adjoint linearized operators.** In order to calculate the minimal norm solution (41), we will need practically useful expressions for the adjoints of the linearized operators of section (4.3). We will present such expressions in this section. The calculation of the actions of these operators will typically require us to solve an adjoint transport problem. This explains the name 'adjoint-field method' of the inversion method employed here.

To start with, a simple calculation gives us the following theorem.

**Theorem 4.1.** *The adjoint operator $\Lambda'[\phi]^*$ which corresponds to the linearized operator $\Lambda'[\phi]$ is given by*

$$\Lambda'[\phi]^* : F \longrightarrow \Phi \quad , \quad \left(\Lambda'[\phi]^* \delta a_{\rm s}\right)(x) = -[\hat{a} - a_{\rm b}(x)] \frac{\delta a_{\rm s}(x)}{|\nabla\phi(x)|} C_\rho(\Gamma)\chi_{B_\rho(\Gamma)}(x). \tag{42}$$

The next theorem describes the adjoint operator $R'_j[a_{\rm s}]^*$ which corresponds to $R'_j[a_{\rm s}]$. Its proof can be found in [9, 10], and is therefore omitted here.

**Theorem 4.2.** *Let $z \in U^*$ be a solution of the (adjoint) transport equation*

$$-\frac{\partial z}{\partial t} - \theta \cdot \nabla z(x,\theta,t) + (a(x) + b(x))\, z(x,\theta,t) - b(x) \int_{S^{n-1}} \eta(\theta \cdot \theta')\, z(x,\theta',t) d\theta'$$

$$= 0 \quad in \quad \Omega \times S^{n-1} \times [0,T] \tag{43}$$

*with 'initial' condition*

$$z(x,\theta,T) = 0 \quad on \quad \Omega \times S^{n-1} \tag{44}$$

*and boundary condition*

$$z|_{\Gamma_+} = \tilde{z}(x,t) := \zeta(x,t). \tag{45}$$

*Then the adjoint operator $R'_j[a_{\rm s}]^*$ is given by*

$$(R'_j[a_{\rm s}]^*\zeta)(x) = I_j(a_{\rm s}) := \int_{[0,T]}\int_{S^{n-1}} u_j(x,\theta,t)z_j(x,\theta,t)\,d\theta dt, \tag{46}$$

*with $u_j$ being a solution of (1)–(3) with source $q_j$ and with $a = a_{\rm b} + a_{\rm s}$. Notice that $\zeta(x,t)$ is applied uniformly into all directions $\theta$ with $\nu \cdot \theta > 0$.* □

Finally, by combining theorems 4.1 and 4.2, we get an expression for the adjoint operator $T'_j[\phi]^*$ which corresponds to the linearized forward operator $T'_j[\phi]$. It is described in the following theorem.

**Theorem 4.3.** *Let $\zeta \in Z$. The adjoint operator $T'_j[\phi]^*$ acts on $\zeta$ in the following way*

$$T'_j[\phi]^*\zeta = \Lambda'[\phi]^*R'_j[\Lambda(\phi)]^*\zeta = \frac{[\hat{a} - a_{\rm b}(x)]}{|\nabla\phi(x)|}C_\rho(\Gamma)\chi_{B_\rho(\Gamma)}(x)I_j(\Lambda(\phi)) \tag{47}$$

*where $I_j(\Lambda(\phi))$ is given by (46) and where $u_j$ solves (1)-(3) and $z_j$ solves (43)–(45) with $a_{\rm s}$ replaced by $\Lambda(\phi)$.*

4.6. **Updating the level set function.** Let us consider the operator $C_{j,n} := T'_j[\phi^{(n)}]T'_j[\phi^{(n)}]^*$ in (41) more closely. In order to calculate a correction $\delta\phi^{(n)}_{\rm MN}$ by (41) we first have to apply the operator $C^{-1}_{j,n}$ to the expression $M_j u_j - \tilde{G}_j$. This operator is difficult to calculate, since typically a large number of forward problems has to be solved for this purpose if the background medium is inhomogeneous and the corresponding Green functions are not easily available. Moreover, the inversion of this operator is usually unstable due to the ill-posedness of the underlying inverse problem. Therefore, in our numerical experiments, we will usually just replace this operator by a suitably chosen approximate operator $\hat{C}_j$ which is easier to calculate and whose inversion is more stable. One possibility for such an approximation (which is adopted here) is to simply use a multiple of the identity, combined with an appropriate selection of receivers in the current step of the nonlinear Kaczmarz-type approach. We refer here for a more thorough discussion of this operator and its approximations to [9, 10, 11, 12, 13, 34]. With this approximation, we define

$$\zeta_j := \hat{C}_j^{-1}\left(M_j u_j - \tilde{G}_j\right). \tag{48}$$

Next, in order to find $\delta\phi^{(n)}_{\rm MN}$, we have to calculate $T'_j[\phi^{(n)}]^*\zeta_j$. Using (47) we get

$$T'_j[\phi^{(n)}]^*\zeta_j = \frac{(\hat{a} - a_{\rm b}(x))}{|\nabla\phi^{(n)}(x)|}C_\rho(\Gamma)\chi_{B_\rho(\Gamma)}(x)I_j(\Lambda(\phi^{(n)})), \tag{49}$$

where $u_j$ solves (1)-(3) with source $q_j$ and $z_j$ solves (43)–(45) with $\zeta = \zeta_j$ and $a_{\rm s} = \Lambda(\phi^{(n)})$.

In the following we want to introduce an approximation to the update formula (49) which in our numerical experiments so far has shown to simplify and stabilize the reconstruction process. We mentioned already earlier that the level set representation of a given contour is not uniquely defined. Therefore, it would be desirable to pick from all those possible level set functions representing a given contour that one with the most convenient properties.

Since we are updating the level set function in each step of our algorithm, it might happen that the denominator in (49) becomes very small, and numerical noise and roundoff-errors can destabilize the inversion routine. In order to avoid this problem, we rescale the level set function $\phi$ after each update with the aim to keep $|\nabla\phi(x)|$ near the boundaries of the inclusions in a certain range of values. As a simple scaling rule, we multiply the level set function with a scaling factor such that the global minimum (or alternatively maximum) of the level set function is always equal to a predefined constant value. Notice that this rescaling procedure does not change the domains which are represented by the given level set function. In other words, we try to pick from the family of level set functions representing a given domain that one which has a convenient behavior of $|\nabla\phi(x)|$ in the neighborhood of the zero level set.

We go one step further and make the assumption that, with this rescaling procedure, the gradient $|\nabla\phi(x)|$ along the boundary can be approximated in (49) by some constant $c_1$. Although, with this approximation, we are not marching anymore into the 'optimal' direction when calculating the updates, we will gain some stability in the inversion routine since we avoid dividing by small numerical derivatives (these are still possible even after the rescaling) in (49) which are strongly influenced by unavoidable numerical noise and roundoff errors. Moreover, our numerical experiments show that these small deviations from the 'optimal' marching direction are easily corrected by the succeeding updates.

With this modification, we get the following *update formula for the level set function*

$$\delta\hat{\phi}^{(n)}(x) = -(\hat{a} - a_{\mathrm{b}}(x))c_1^{-1}\, C_\rho(\Gamma)\chi_{B_\rho(\Gamma)}(x)I_j(\Lambda(\phi^{(n)})). \qquad (50)$$

In (50) we did not incorporate the rescaling procedure described above.

Notice that, although we did not explicitly impose any regularity constraints on the updates (50), they are now in the range of $R'_j[\Lambda(\phi^{(n)})]^*$ (up to the factor $\hat{a} - a_{\mathrm{b}}(x)$) which implicitly gives us some information about the regularity we can expect. Our numerical experiments so far indicate that the degree of regularity which is achieved by applying (50) is typically sufficient 'for practical purposes' in our numerical examples.

4.7. **Sensitivity analysis.** A valuable tool for designing an imaging experiment is a careful sensitivity analysis of the underlying inverse problem. The sensitivity functions which correspond to a nonlinear inverse problem are functions which decompose the linearized residual operator (which is often called 'Fréchet derivative') as well as its adjoint operator, as described in more details in [14, 15]. The sensitivity functions corresponding to one given source position and one given receiver location describe the sensitivity of this source-receiver pair to small variations of the parameters in the medium. Locations in the medium where the sensitivities are small have almost no influence on the measurements at this specific detector location. Therefore, hidden objects at these locations are difficult to recover unless other source-receiver pairs have larger sensitivity values at these specific locations. Ideally, we would like to design an experiment where certain source-receiver configurations have a sensitivity structure which is large on a very small area in the domain, and very small elsewhere. The corresponding data can then be used for determining the coefficients at these specific locations. If we find a set of sensitivity functions such that each location in the medium is 'probed' by such an ideal configuration, we expect to get good and fast reconstructions. Unfortunately, in diffusion-type inverse problems like ours, the sensitivity functions are typically spread out over a relatively large area. In addition, they have a large dynamical range since they have very large values close to the source and receiver locations, and very small values in the regions of interest. Therefore, it is very difficult to probe localized areas inside the domain with suitable superpositions of sensitivity functions. The role of the sensitivity functions in our TBT-reconstruction scheme is described in the following theorem.

**Theorem 4.4. Transport sensitivity functions.** *Given a source (6). We can find functions $\Psi_a(x_r, t_r; x_{sc})$ with the following property.*
*1.)'Projection step': The action of the linearized residual operator $R'[a_s]$ on the perturbation $\delta a$ in parameter space can be described as follows*

$$\left(R'[a_s]\delta a\right)(x_r, t_r) = \int_\Omega \Psi_a(x_r, t_r; x_{sc})\delta a(x_{sc})dx_{sc}. \tag{51}$$

*2.) 'Backprojection step': The action of the adjoint linearized residual operator $R'[a_s]^*$ on the vector $\zeta$ in data space can be described as follows*

$$\left(R'^*_a \zeta\right)(x_{sc}) = \int_{\partial\Omega} \int_{[0,T]} \Psi_a(x_r, t_r; x_{sc})\zeta(x_r, t_r)\, dt_r dx_r. \tag{52}$$

*The functions $\Psi_a(x_r, t_r; x_{sc})$ are called 'absorption transport sensitivity functions'. Here, $x_r$ and $t_r$ denote receiver location and receiver time, respectively, and $x_{sc}$ is an arbitrary scattering point in the medium.*

A derivation of this theorem, and a more detailed discussion of these sensitivity functions, can be found in [14]. We have displayed typical sensitivity functions for our imaging configuration in figure 1. The upper row of the figure shows the distribution of the scattering parameter for two different imaging situations. On the left, we have no clear layers, on the right we have a clear layer close to the boundary of the imaging domain which might represent for example the Cerebro-Spinal Fluid Layer (CFL) of the human head [1, 35, 39]. The source position considered here is indicated as point 'P1' in the figure on the left hand side. The receiver is located at the point 'P3'. (We are not considering here the also marked point 'P2'.) These points are the same in the situation shown on the right hand side. The second row displays the sensitivity functions with respect to the absorption coefficient for these two situations and for those data which are due to a source position at P1 and taken at the receiver position P3 at an early time, namely $t = 10$ s (recall the normalization $c = 1$ cm/s). It can be seen that these early times are very sensitive to the clear layer. Our interpretation of this observation is that the early photons are very likely to have been channeled through this clear layer, and are therefore very sensitive to small variations of the absorption parameter inside the clear layer.

The third row shows the same sensitivity functions, but at a later time $t = 24$ s. We see that these later times are more sensitive to variations in the absorption parameter inside the domain and less to those inside the clear layer. Photons which arrive that late at the given receiver, are very likely to have travelled a long time in the scattering medium, and are therefore more affected by parameter variations inside the medium. However, we also observe that these later sensitivity functions are more 'spread-out' than those corresponding to earlier times, since the photons had more time to reach remote areas inside the imaging domain.

Table 1: *The level set reconstruction algorithm.*

*Initialization.*

$n = 0$;

$(D^{(0)}, a^{(0)}, \phi^{(0)})$ given from TBT.

*Reconstruction loop.*

FOR $i = 1 : I_m$     perform $I_m$ sweeps

  FOR $j = 1 : p$     march over source positions

  $\zeta_j = \hat{C}_j^{-1}(M_j u_j - \tilde{G}_j);$     $u_j$ solves (1)–(3) with $a^{(n)}$, $q_j$

  $\delta\phi^{(n)} = -(\hat{a} - a_b(x))I_j(\Lambda(\phi^{(n)}))\chi_{B_\rho(\Gamma)};$     $z_j$ solves (43)-(45) with $a^{(n)}$ and $\zeta_j$

  $\phi^{(n+1)} = C_{LS}^{(n)}(\phi^{(n)} + \eta\frac{C_\rho(\Gamma)}{c_1}\delta\phi^{(n)});$     update level set function

  $a^{(n+1)} = a_b + \Lambda(\phi^{(n+1)});$   $n = n+1;$     reinitialization $n \to n+1$

  END

END

$(D^{(N)}, a^{(N)}, \phi^{(N)}) = (D[\phi^{(n)}], a_b + \Lambda(\phi^{(n)}), \phi^{(n)});$     final reconstruction.

Since the sensitivity structure of these later photons is more focused onto the domain of interest, we will use only these data for solving our inversion task. The inversion process incorporates an iterative backprojection ('backtransport') step along these sensitivity functions, as indicated in theorem 4.4 and explained in [14, 15]. Since the functions along which the data are backprojected are very spread-out, we expect a loss in resolution compared to those reconstructions which can make use of the more focused earlier-time sensitivity functions.

4.8. **The algorithm.** The nonlinear Kaczmarz-type method for shape reconstruction using level sets is described in brief algorithmic form in table 1. Here, $\eta$ is a relaxation parameter for the update of the level set function which is determined empirically. The operator $\hat{C}_j^{-1}$ is taken to be simply a multiple of the identity, such that it can be considered as part of $\eta$. The constant $C_\rho(\Gamma)$ could be calculated explicitly for the actual curve $\Gamma^{(n)}$, or it could be approximated by some value corresponding to a simple geometrical object (to give an example, in case of a single circle it would be $C_\rho(\Gamma) = (2\rho)^{-1}$). In our numerical experiments so far, however, it is simply considered as part of $\eta$. The same holds true for $c_1$. The constant $\rho$ is in our numerical experiments chosen between one and two grid cells. The scaling factor $C_{LS}^{(n)}$ is determined after each update to keep the global minimum (or maximum) of the level set function at a constant value.

Notice that we typically apply a very small relaxation parameter $\eta$, which has the effect that the accumulated contributions of the earlier updates are not totally destroyed by the new correction, but instead have some 'inertia' which is the larger the smaller we choose $\eta$. This slows down the whole reconstruction process to some degree, but, on the other hand, stabilizes it significantly. In our numerical experiments, we typically choose $\eta$ such that the corrections to the boundaries of the shapes in each update are not much larger than one or two pixels in the normal direction to these boundaries.

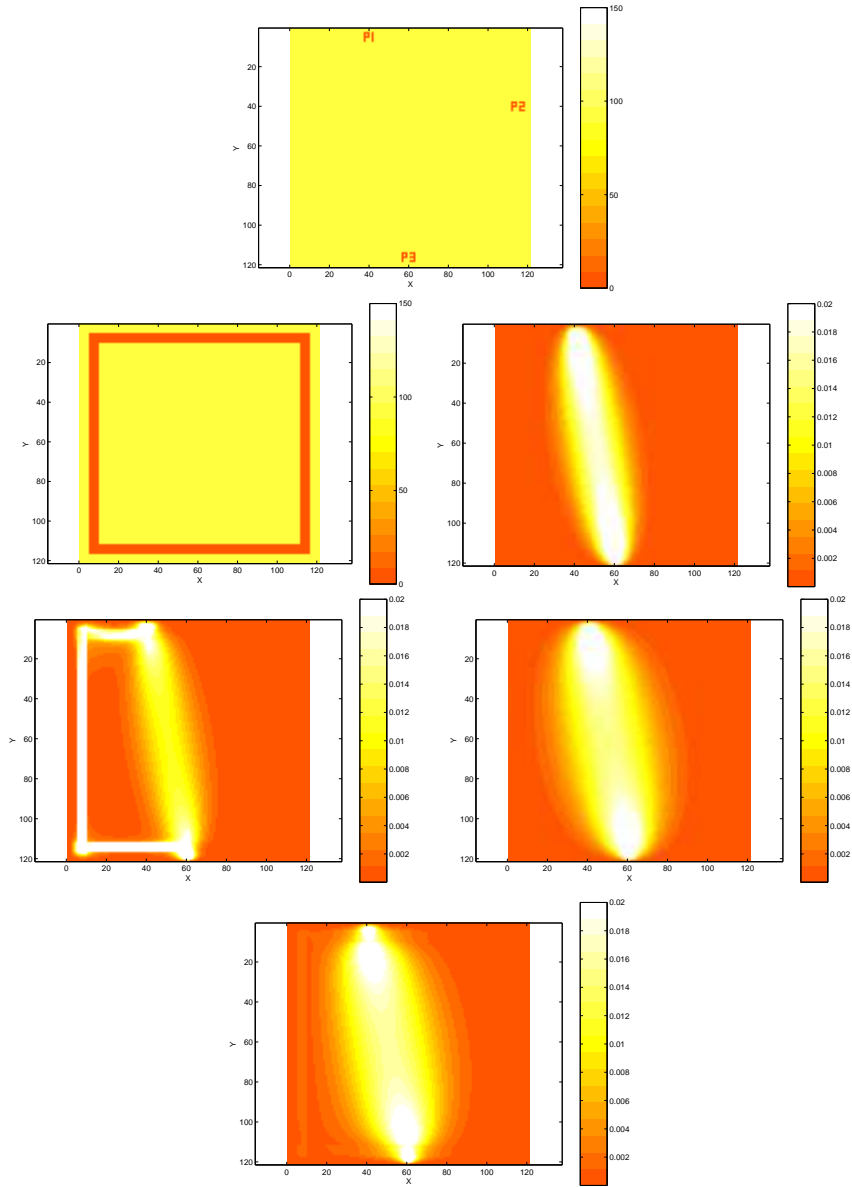

FIGURE 1.    Transport absorption sensitivity functions $-\Psi_a(x_r, t_r; x_{sc})$ for a source located at $P1$ and a receiver at $P3$. Left column: homogeneous background. Right column: homogeneous background with clear region close to boundary. Upper row: the domains. Middle row: at an early time step ($t_r = 10$ s). Bottom row: at a later time step ($t_r = 24$ s).

## 5. THE TBT-METHOD AND THE FIRST GUESS

For starting our shape reconstruction method using level sets we will need an initial guess $(D^{(0)}, a^{(0)}, \phi^{(0)})$. Although it is possible just to create an arbitrary initial guess without using any data at all, there are several reasons for employing a preprocessing step which is relatively inexpensive and which yields a good initial guess for the level set reconstruction method. This preprocessing routine usually indicates already the approximate locations of possible inclusions, such that the amount of work which has to be done by the level set method is significantly reduced, and the algorithm converges much faster.

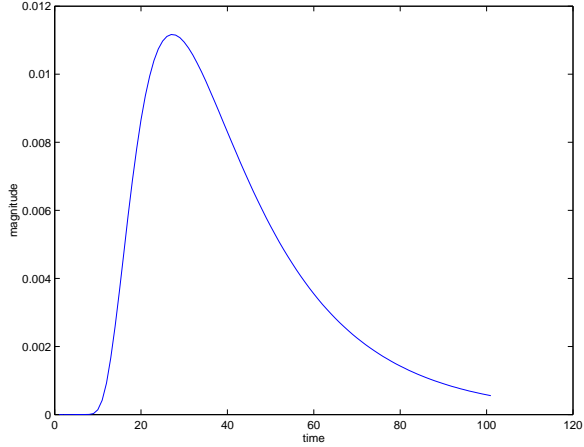

Figure 2. A typical data curve at one receiver in DOT.

We are employing an inverse scattering algorithm (the TBT-algorithm) for this purpose which was introduced in [9, 10] and which shares basic features with the level set reconstruction scheme. Typically, in the diffusive regime, this inverse scattering algorithm yields already in the early sweeps good low-contrast approximations to the high-contrast objects. These low-contrast approximations usually indicate already the approximate locations of the unknown objects, but typically do not capture the correct number, sizes, shapes or contrasts of these objects. For this, many more iterations of the TBT scheme would be necessary. Therefore, it is advantageous to apply only a few steps of this inverse scattering method, and to switch then to a shape reconstruction formulation which explicitly takes into account the high-contrast assumption. We will show that in fact much better reconstructions can be achieved rapidly in these situations by the shape reconstruction method compared to the inverse scattering algorithm. We refer for a detailed description of the TBT-method to the publications [9, 10], and describe in the following only the scheme which we use for extracting a suitable initial triple $(D^{(0)}, a^{(0)}, \phi^{(0)})$ from the TBT reconstructions $a_{\text{TBT}}(x)$.

Let as before $\tilde{D}$ denote the sought domain. We assume that we are working in a high-contrast situation, such that exactly one of the following conditions is satisfied

$$\hat{a} - a_{\text{b}}(x) \gg 0 \quad \text{for all } x \in \tilde{D} \tag{53}$$

$$\hat{a} - a_{\text{b}}(x) \ll 0 \quad \text{for all } x \in \tilde{D}. \tag{54}$$

Since we know $\hat{a}$ and $a_{\text{b}}(x)$, we know the constant

$$\text{sign}(\tilde{D}) := \begin{cases} 1 & , \text{ if (53) holds} \\ -1 & , \text{ if (54) holds}. \end{cases} \tag{55}$$

Choose a threshold value $0 < \gamma_{\text{LS}} < 1$ (in our numerical examples presented in section 6 we use $\gamma_{\text{LS}} = 0.9$) and define

$$a_{\text{LS}} := \gamma_{\text{LS}} \max_{x \in \Omega} a_{\text{TBT}}(x). \tag{56}$$

For the level set zero $L_0^{(0)}$ of $\phi^{(0)}$ we require that

$$L_0^{(0)} = \{x \in \Omega : \quad a_{\text{TBT}}(x) = a_{\text{LS}}\}. \tag{57}$$

This means that we want all points of $\Omega$ where the reconstruction $a_{\text{TBT}}(x)$ has exactly the value $a_{\text{LS}}$ to be mapped to zero by the level set function $\phi^{(0)}$

$$\phi^{(0)}(x) = 0 \quad \text{for all } x \in L_0^{(0)}. \tag{58}$$

The level set function is now defined as

$$\phi^{(0)}(x) = C_{\text{LS}}^{(0)} \text{sign}(\tilde{D}) \left(a_{\text{LS}} - a_{\text{TBT}}(x)\right), \tag{59}$$

where $C_{\text{LS}}^{(0)}$ is some suitably chosen scaling factor. Notice that (58) and (59) are consistent. The initial scattering domain $D^{(0)}$ and the absorption $a^{(0)}$ are defined as

$$D^{(0)} = D[\phi^{(0)}], \qquad a^{(0)} = a_{\text{b}} + \Lambda(\phi^{(0)}). \tag{60}$$

Together with $\phi^{(0)}$ they form an admissible triple $(D^{(0)}, a^{(0)}, \phi^{(0)})$.

## 6. Numerical Experiments

**6.1. The computational domain.** In our numerical experiments we consider a domain of the size $5 \times 5$ cm$^2$ which is divided into $50 \times 50$ pixels each having a size of $1 \times 1$ mm$^2$. One time step lasts 0.2 seconds ([s]), and measurements are taken during the first 100 time steps or 20.0 s. Since the velocity of the light in the medium is normalized to $c = 1.0$ cm/s, each particle can travel 20 cm in that time. The directional variable is discretized into 12 direction vectors which are equidistantly distributed over the unit circle. The scattering law is a Henyey-Greenstein function (5) with $g = 0.9$. The space and time derivatives of (1) are discretized by using a finite difference scheme. A more detailed description of our numerical scheme for solving (1)–(3) can be found in [9, 10].

Our sources have a width of 5 pixels or 5 mm along the boundary. Using these slightly extended sources reduces artifacts which tend to appear close to the boundaries, but is not essential for the presented method. In total, there are 16 of these sources, 4 at each side of the domain. These sources are emitting an ultrashort (one time step) laser pulse at time $t = 0$ in the direction perpendicular to the boundary into the medium. They are arranged in a way such that only the central 20 mm of each side are illuminated during the whole imaging process. In other words, those parts of the boundary whose distance to any of the corners is less than 15 mm are not illuminated in our experiments. The reason for this source arrangement is the location of the clear layer close to the boundary, combined with the fact that our domain is rectangular rather than circular. By this source arrangement, we avoid that ballistic (i.e. unscattered) photons travel along the clear regions.

Receivers are positioned at those parts of the boundary which have a distance (along the boundary) of at least 5 cm or 50 pixels to the source position. Therefore, only half of the boundary is actually covered by the receivers in a given experiment.

This arrangement has been chosen since the source-near receivers are not very sensitive to the regions of interest inside the domain. Therefore, they are not supposed to carry much information about the sought objects. In addition, we believe that the approximation of the matrices $C_{j,n}$ by a multiple of the identity is easier to justify with this source-receiver geometry, because of the reduced 'dynamic range' of the corresponding data. However, we do not have any proof for this assumption so far.

Another important detail is that we are only using those photons as data which arrived at the receiver positions at times between $t = 8$ s and $t = 20$ s. In other words, the early photons are not used for the inversion task, since they are more likely to have been channeled through the clear region close to the boundary. We plan to use the information of these early photons in our future work for the additional determination of the structure of the clear layers from the measured data, which is here simply assumed to be known.

The background of our imaging medium is homogeneous except of the embedded clear layers. These homogeneous regions have a scattering parameter of $b = 100.0$ cm$^{-1}$, and an absorption parameter of $a_b = 0.1$ cm$^{-1}$. Inside the clear regions, the parameters are $b = 0.01$ cm$^{-1}$ and $a_b = 0.01$ cm$^{-1}$. In the first two numerical experiments, we have only one clear layer close to the boundary and surrounding the interior imaging domain, see figures 1, 3, 4. In the third numerical experiment, we have two additional disc-shaped clear regions imbedded in the interior imaging domain, see figure 5.

6.2. **Three numerical experiments.** The results of our first numerical experiment are shown in figure 3. We want to recover three disc-shaped objects in the interior of the domain from time-dependent data. The true configuration is shown in the bottom right image of the figure. The values inside the obstacles are $a = 0.5$ cm$^{-1}$ and $b = 100.0$ cm$^{-1}$. The data are created by running our finite differences discretization of the linear transport equation (1)-(3) on the correct model. So far, no noise is added to the data. We plan to investigate the effect of noise in the data more carefully in our future research.

First, we run our preprocessing step (a TBT algorithm) in order to find a good estimate of the number and locations of the obstacles. The results of these runs after 5 sweeps and after 20 sweeps of the TBT method are shown in the first two images of the figure. Notice that one sweep of the TBT method consists of using the data for each of the sources exactly once, in a nonlinear Kaczmarz-type scheme. For more details of that scheme see section 4.4 and [9, 10, 34]. We see that the locations and even the number of the unknown obstacles are already indicated in the reconstructions after 20 sweeps, but that the contrast and the shapes are far away from being correct. We use the result in order to start our level-set based shape reconstruction scheme. The figure shows the initial shape which is extracted from the TBT reconstruction, the reconstruction after six steps of the algorithm (which means after only using the information corresponding to six of the 16 source positions once), and after 16 steps or one sweep of the algorithm. We see that after only one additional sweep, now with the level set based shape reconstruction scheme, we have arrived at a very good estimate of the shapes and locations of the unknown obstacles. Additional iterations (sweeps) do not improve the reconstruction significantly, such that they are not displayed here. In this second step we have assumed that we know the correct contrast of the obstacles to the background. We conclude that the explicit incorporation of this additional information of the contrast in our reconstruction scheme enables us to get rapidly a much better estimate of the shapes of the obstacles than it would have been possible by just continuing the pixel-based TBT reconstruction scheme.

Certainly, the final goal of the reconstruction algorithm will be to find the shapes as well as the contrast of the obstacles simultaneously from the given data. We will address this question in our future research. Notice that the shapes of the obstacles change topology during the evolution process. The level set method has no difficulties in keeping track of these topology changes.

As mentioned, in practical applications we usually do not know the correct values of the absorption parameter inside the shapes. In our second numerical experiment, we investigate the behavior of the level set scheme in the case that we prescribe a slightly incorrect value for the contrast. Instead of $a_s = 0.5$ cm$^{-1}$ inside the obstacle we assume $a_s = 0.55$ cm$^{-1}$. The new configuration is displayed in the bottom right image of figure 4. The other images of the figure show the results of the algorithm. We see from these results that the effect of this model error on the reconstruction is very small. After only one additional sweep with the level set based reconstruction scheme, we arrive at a very good estimate for the shapes of the obstacles. These shapes also do not differ much from those where the correct value is used. That indicates that the data which are used here are not very sensitive to small changes in the parameter values.

In our third numerical experiment we investigate a slightly more complex situation. Now, there are two additional disc-shaped clear regions in the domain as shown in Figure 5. The bottom right image of the figure shows the true configuration. Although we assume that we know the boundaries of these additional clear regions, such that we do not have to recover them simultaneously with the unknown obstacles, these regions might have a negative impact on our ability to recover these obstacles due to the changed sensitivity structure. We refer to [14] for visualizations of the sensitivity functions in a similar situation. As an additional difficulty, we assume now that each of the three objects which we are looking for has a different contrast to the background. The three values of the absorption parameter inside the objects are $a = 0.4 \,\text{cm}^{-1}$, $a = 0.5 \,\text{cm}^{-1}$, and $a = 0.6 \,\text{cm}^{-1}$. In our reconstruction process, we apply the constant value $a = 0.4 \,\text{cm}^{-1}$ for all three unknown obstacles, such that two of the objects are reconstructed with an incorrect value. Figure 5 shows first the reconstruction of the TBT-method after 20 sweeps. This is used for extracting a first guess for our shape reconstruction scheme, as it is shown in the second image on the top row. The third image on the top row shows the reconstruction after 6 steps of the shape reconstruction scheme (which means that only the information corresponding to the first six source position has been used exactly once), after three sweeps, and after 10 sweeps of the scheme. As in our first two numerical examples, the algorithm delivered very good estimates of the shapes after only a few sweeps. The evolution of the norm of the residuals during the reconstruction is displayed in figure 6. The residuals converge now more slowly than in the first two examples, and the later sweeps still yield small improvements to the reconstructions. Notice also that the final reconstructions of the obstacles tend to be slightly more extended than the original objects in the case that a smaller contrast value was used for the inversion. This, of course, might have been expected since a larger area of absorption material is now necessary in order to have the same effect on the data as the original object with a higher absorption value. Finally, we observe also here that the shapes do change topology in order to arrive at the final solution. This is handled easily by the level set formulation.

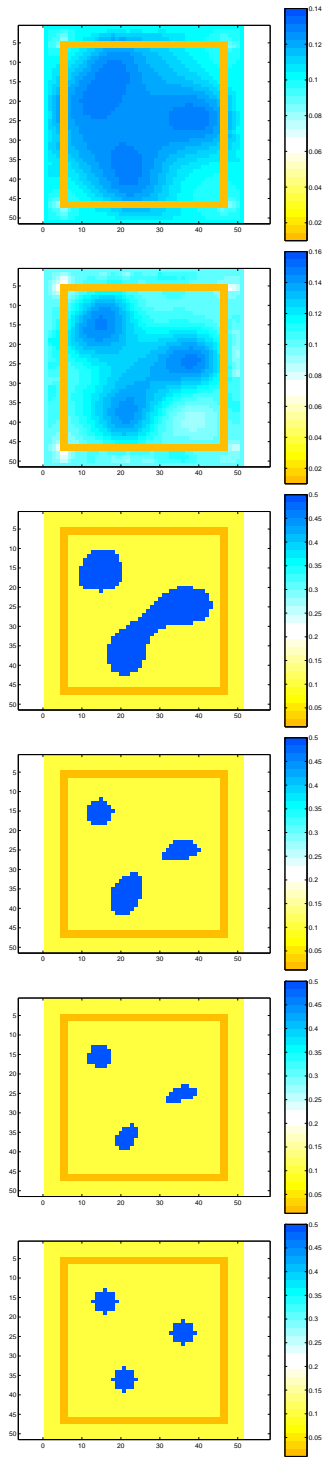

FIGURE 3. The correct value $a = 0.5\,\text{cm}^{-1}$ is used for the reconstruction. Top row from left to right: $a_{\text{TBT}}$ after 5 sweeps of TBT; after 20 sweeps of TBT; initial shape for level set algorithm. Bottom row from left to right: shape after 6 steps of level set algorithm; after 16 steps (i.e. 1 sweep); true object.

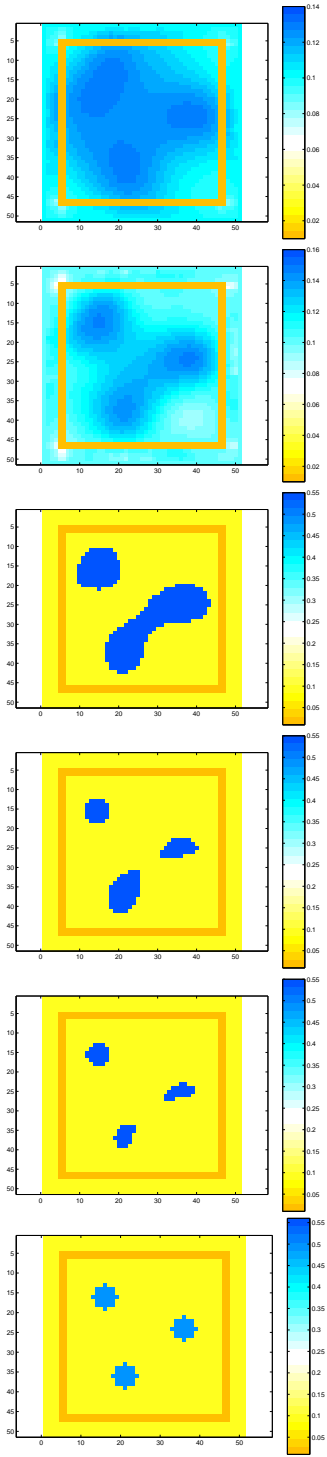

FIGURE 4. An approximate value $a = 0.55\,\text{cm}^{-1}$ (instead of $a = 0.5\,\text{cm}^{-1}$) is used for the reconstruction. Top row from left to right: $a_{\text{TBT}}$ after 5 sweeps of TBT; after 20 sweeps of TBT; initial shape for level set algorithm. Bottom row from left to right: shape after 6 steps of level set algorithm; after 16 steps (i.e. 1 sweep); true object.

## 7. Summary and future research directions

We have presented a novel two-step shape reconstruction method for DOT which uses level sets and adjoint fields. The first step of the method employs the TBT method in order to determine a first guess for the shapes. Then, the second step uses this first guess to start a shape evolution scheme in order to complete the inversion task. The combination of these two schemes makes use of the fact that the TBT method typically finds a good guess for the locations of hidden obstacles in the early sweeps, and that the level set based shape reconstruction method is particularly fast once a good approximation for the shapes has been found. Therefore, the combination of these two methods enables us to get very fast and accurate reconstructions of the unknown objects. The reformulation of the inverse problem as a shape reconstruction problem can be considered as a regularization technique. The incorporated prior information about the high contrast of the inclusions and the approximately known values of the parameters stabilizes and speeds up the inversion task. We believe that also the simultaneous reconstruction of shapes and contrast will be possible with a slight modification of the algorithm, and plan to investigate this interesting problem in our future work.

So far, we have been using mainly the information from later arriving photons for the inversion task, in order to cope with the specific structure of the sensitivity functions when clear layers are present. In our future work, the additional information of the early photons will be used in order to try to determine also the structure of the clear layers from the given data. Moreover, we plan to extend the imaging problem to finding the scattering parameter $b(x)$ as well. An important question to be addressed is also the proper incorporation of additional regularization techniques, and the correct choice of the underlying function spaces. Especially the choice of the level set space is not obvious from the beginning. The proper choice of this space will ideally also be connected to our regularization strategy. We conclude that the novel two-step shape reconstruction technique presented here is extremely promising for the application in DOT, and that it poses many interesting and mathematically challenging questions which are open for future research.

## References

[1] Arridge S R (1999) Optical tomography in medical imaging *Inverse Problems* **15** (2), R41–R93
[2] Burger M (2001) A level set method for inverse problems *Inverse Problems* **17** 1327–1355
[3] Case K M and Zweifel P F (1967) *Linear Transport Theory* (New York: Plenum Press)
[4] Claerbout J F 1976 *Fundamentals of Geophysical Data Processing: With Applications to Petroleum Prospecting* McGraw-Hill, New York
[5] Dautray R and Lions J L (1993) *Mathematical Analysis and Numerical Methods for Science and Technology* Vol.6 (Berlin: Springer)
[6] Dierkes T (2000) Rekonstruktionsverfahren zur optischen Tomographie *Thesis* Preprints Angewandte Mathematik und Informatik 16/00-N, Münster (in german)
[7] Dierkes T, Dorn O, Natterer F, Palamodov V and Sielschott H (2001) Fréchet derivatives for some bilinear inverse problems *to appear in 'SIAM Journal on Applied Mathematics'*.
[8] Dines K A and Little R J 1979 Computerized geophysical tomography *Proc IEEE* **67** 1065–73
[9] Dorn O (1997) Das inverse Transportproblem in der Lasertomographie *Thesis* Preprints "Angewandte Mathematik und Informatik" 7/97-N, Münster (in german)

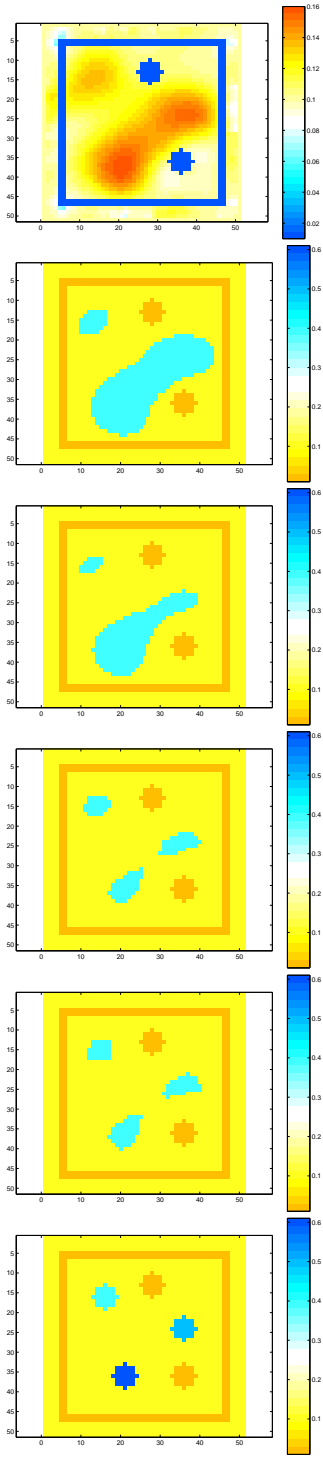

FIGURE 5. The smallest ($a = 0.4\,\mathrm{cm}^{-1}$) of the three values ($a = 0.4\,\mathrm{cm}^{-1}$, $a = 0.5\,\mathrm{cm}^{-1}$ and $a = 0.6\,\mathrm{cm}^{-1}$) is used for the reconstruction. Top row from left to right: $a_{\mathrm{TBT}}$ after 20 sweeps of TBT; initial shape for level set algorithm; shape after 6 steps of level set algorithm; Bottom row from left to right: shape after 3 sweeps of level set algorithm; after 10 sweeps; true object.

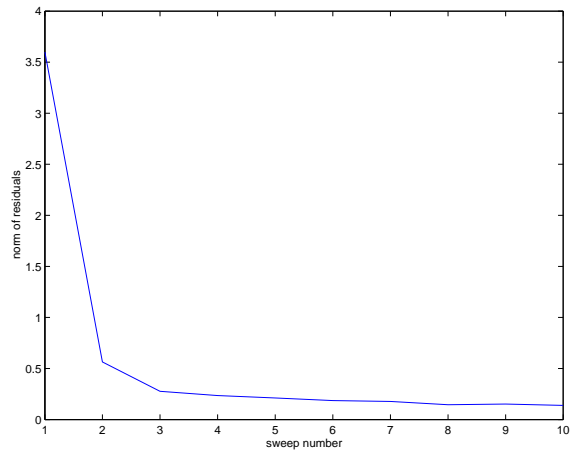

Figure 6. Norm of the residuals during the level set reconstruction process shown in figure 5

[10] Dorn O (1998) A transport-backtransport method for optical tomography *Inverse Problems* **14** 1107–1130
[11] Dorn O, Bertete-Aguirre H, Berryman J G and Papanicolaou G C 1999 A nonlinear inversion method for 3D electromagnetic imaging using adjoint fields *Inverse Problems* **15** 1523-1558
[12] Dorn O, Miller E L and Rappaport C M 2000 A shape reconstruction method for electromagnetic tomography using adjoint fields and level sets, *Inverse Problems* **16** 1119-1156
[13] Dorn O, Miller E L and Rappaport C M 2001 Shape reconstruction in 2D from limited-view multifrequency electromagnetic data, AMS series Contemporary Mathematics, **278**, 'Radon Transform and Tomography', 97-122
[14] Dorn O 2000 Scattering and absorption transport sensitivity functions for optical tomography *Optics Express* **7** (13) (Focus Issue: Diffuse Optical Tomography) 492–506, electronically available at http://www.opticsexpress.org/opticsexpress/tocv7n13.htm
[15] Dorn O, Bertete-Aguirre H, Berryman J G and Papanicolaou G C 2001 Sensitivity analysis of a nonlinear inversion method for 3D electromagnetic imaging in anisotropic media *to appear in 'Inverse Problems'*.
[16] Elangovan V and Whitaker R (2001) From Sinograms to surfaces: a direct approach to the segmentation of tomographic data *International Conference on Medical Image Computing and Computer-Assisted Intervention (MICCAI) 2001* Utrecht, The Netherlands
[17] Gaudette J J, Brooks C A, Kilmer M E, Miller E L, Gaudette T and Boas D (2000) A comparison study of linear reconstruction techniques for diffuse optical tomographic imaging of absorption coefficient, *Phys.Med.Biol.* **45** (4) 1051–1070.
[18] Grünbaum F A (2001) A nonlinear inverse problem inspired by three-dimensional diffuse tomography, *Inverse Problems* **17** (6) 1907-1922
[19] Harabetian E and Osher S (1998) Regularization of ill-posed problems via the level set approach *SIAM J.Appl.Math.* **58** (6) 1689–1706
[20] Hielscher A, Klose A and Beuthan J Evolution strategies for optical tomographic characterization of homogeneous media, *Optics Express* **7** (13) (Focus Issue: Diffuse Optical Tomography) 507
[21] Hintermüller M and Ring W (2001) A level set approach for the solution of a state-constraint optimal control problem *Report No. 212* University of Graz and Technical University of Graz, February 2001
[22] Ito K, Kunisch K and Li Z (2001) Level-set function approach to an inverse interface problem *Inverse Problems* **17** 1225–1242
[23] Jackowska-Strumillo L, Sokolowski J, Zochowski A and Henrot A (2001) On numerical solution of shape inverse problems, preprint.

[24] Kaltenbach J M and Kaschke M (1993) Frequency- and Time-Domain Modelling of Light transport in Random Media, in: *Medical Optical Tomography* ed. Potter (1993) SPIE Optical Engineering Press, Vol. IS11 65-86
[25] Klibanov M V, Lucas T R and Frank R M (1997) A fast and accurate imaging algorithm in optical/diffusion tomography *Inverse Problems* **13** 1341-1361
[26] Kolehmainen V, Arridge S R, Lionheart W R B, Vauhkonen M and Kaipio J P 1999 Recovery of region boundaries of piecewise constant coefficients of an elliptic PDE from boundary data *Inverse Problems* **15** 1375-1391
[27] Kolehmainen V (2001) Novel approaches to image reconstruction in diffusion tomography *thesis*, Kuopio University Publications C. Natural and Environmental Sciences 125
[28] Litman A, Lesselier D and Santosa F 1998 Reconstruction of a two-dimensional binary obstacle by controlled evolution of a level-set *Inverse Problems* **14** 685-706
[29] Miller E, Kilmer M and Rappaport C 1999 A New Shape-Based Method for Object Localization and Characterization from Scattered Field Data, submitted to *IEEE Trans. Geoscience and Remote Sensing*
[30] Moscoso M, Keller J B, and Papanicolaou G (2001) Depolarization and blurring of optical images by biological tissues, *J. Opt. Soc. Am. A* **18** 948-960
[31] Natterer F 1986 *The Mathematics of Computerized Tomography* (Stuttgart: Teubner)
[32] Natterer F and Wübbeling F 1995 A propagation-backpropagation method for ultrasound tomography *Inverse Problems* **11** 1225–1232
[33] Natterer F 1996 Numerical Solution of Bilinear Inverse Problems *Preprints* "Angewandte Mathematik und Informatik" 19/96-N Münster. This preprint is electronically available at http://wwwmath.uni-muenster.de/math/num/Preprints/1997/natterer_3.
[34] Natterer F and Wübbeling F, *Mathematical Methods in Image Reconstruction*, Monographs on Mathematical Modeling and Computation 5, SIAM 2001.
[35] E. Okada, M. Firbank, M. Schweiger, S. R. Arridge, M. Cope and D. T. Delpy, "Theoretical and experimental investigation of near-infrared light propagation in a model of the adult head," Appl. Opt. **36** (1), 21–31 (1997)
[36] Osher S and Sethian J 1988 Fronts propagation with curvature dependent speed: Algorithms based on Hamilton-Jacobi formulations *Journal of Computational Physics* **56** 12-49
[37] Osher J S and Santosa F (2001) Level set methods for optimization problems involving geometry and constraints I. Frequencies of a two-density inhomogeneous drum *J. Comp. Phys.* **171** 272-288
[38] Ramananjaona C, Lambert M, Lesselier D and Zolesio J-P (2001) Shape reconstruction of buried obstacles by controlled evolution of a level set: from a min-max formulation to numerical experiments *Inverse Problems* **17** 1087–1111
[39] J. Riley, H. Dehghani, M. Schweiger, S. R. Arridge, J. Ripoll and M. Nieto-Vesperinas (2000) "3D Optical Tomography in the Presence of Void Regions," *Optics Express* **7** (13) (Focus Issue: Diffuse Optical Tomography) 462
[40] Santosa F 1996 A Level-Set Approach for Inverse Problems Involving Obstacles *ESAIM: Control, Optimization and Calculus of Variations* **1** 17-33
[41] Sethian J A 1990 Numerical algorithms for propagating interfaces: Hamilton-Jacobi equations and conservation laws *J. Diff. Geom.* **31** 131-61
[42] Sethian J A 1999 Level Set Methods and Fast Marching Methods (2nd ed) Cambridge University Press
[43] Sethian J A and Wiegmann A (1999) Structural boundary design via level set and immersed interface methods J. Comp. Phys. **163** (2) 489-528
[44] Sielschott H (2001) Rückpropagationsverfahren für die Wellengleichung in bewegtem Medium *thesis*, Preprints "Angewandte Mathematik und Informatik" University of Münster, Germany (in german)
[45] Svaiter B F and Zubelli J P (2001) Convexity of the diffuse tomography model *Inverse Problems* **17** 729–738
[46] Zhao H K, Osher S, Merriman B, and Kang M (2000) Implicit and Non-parametric Shape Reconstruction from Unorganized Points Using Variational Level Set Method *Computer Vision and Image Understanding* **80** 295-319

Oliver Dorn, Department of Computer Science, University of British Columbia, 201-2366 Main Mall, Vancouver, B.C. V6T1Z4, Canada

*E-mail address*: `dorn@cs.ubc.ca`